\documentclass[12pt]{article}

\usepackage{hyperref}

\usepackage{amsfonts, amsmath, amssymb, amsthm}
\usepackage{a4}

\usepackage{blkarray}
\usepackage{bbm}
\usepackage{upgreek}

\usepackage{bm}

\usepackage{graphicx}

 \usepackage{xcolor}

 \definecolor{Refkey}{RGB}{255,127,0}
 \definecolor{Labelkey}{RGB}{127,0,255}
 \makeatletter 
  \def\SK@refcolor{\color{Refkey}}
  \def\SK@labelcolor{\color{Labelkey}}
 \makeatother

  \definecolor{mdg}{RGB}{0,177,0} 
  \definecolor{mdb}{RGB}{0,0,191}
  \definecolor{mddb}{RGB}{0,0,91}
  \definecolor{mdy}{RGB}{255,69,0} 
  \definecolor{gray}{RGB}{99,99,99}

\theoremstyle{definition}

\newtheorem{example}{Example}

\theoremstyle{remark}

\allowdisplaybreaks[3]

\title{Solutions to quantum tetrahedron equation with two colors and nonnegative matrix entries}
\author{Igor G. Korepanov}

\date{April--May 2024}

\begin{document}

 \sloppy

\maketitle

\begin{abstract}
In this short note, we construct solutions to quantum tetrahedron equation of the kind ``with variables on the edges''. Each of these variables takes just two values, called sometimes ``colors''. We propose two different constructions. The first of them involves, in particular, two $\mathcal R$-operators each depending on one parameter, while these parameters are independent from each other; the number of nonvanishing matrix entries in these two $\mathcal R$-operators is either 12 or~14. In the second construction, we use what may be called ``modified tetrahedron in a direct sum'' to produce, again, quantum tetrahedron solutions. All matrix entries of our $\mathcal R$-operators are nonnegative, if the relevant parameters are chosen properly.
\end{abstract}

\section{Introduction}\label{s:i}

Tetrahedron equation is intended for constructing exactly solvable models of three-dimen\-sional statistical physics, as well as two-dimen\-sional quantum lattices. It was invented, in the form ``with variables on the faces'', by Zamolodchikov~\cite{Zamolodchikov}, who gave its first nontrivial solution. Here we deal with the equation ``with variables on the edges''; its first nontrivial solution was found in~\cite{v}.

One problem with tetrahedron equation is that statistical physics requires \emph{positive} Boltzmann weights. This means that the entries of at least one matrix entering in the solution of tetrahedron equation---these are ``local Boltzmann weights''---must be either positive or zero (supposing that zeros could be interpreted as corresponding to some prohibited configurations).

Here we propose some unusual solutions of quantum tetrahedron equation, with nonnegative matrix elements, and where the variable on each edge takes just two values---called sometimes ``colors''. For these colors, we take the elements $0,1 \in \mathbb F_2$ of the smallest Galois field.

Note that a tetrahedron solution with positive Boltzmann weights, but with an \emph{infinite} number of colors, was earlier proposed in~\cite{MBS}.

Below, we identify linear operators with their matrices, as there always is a fixed basis.

\section{Tetrahedron equation in a direct sum and permutation-type solutions to quantum tetrahedron equation}\label{s:d}

\emph{Tetrahedron equation in a direct sum} will mean, in this short note, the following relation between four $3\times 3$ matrices $R^{(1)}$, $R^{(2)}$, $R^{(3)}$ and $R^{(4)}$ with entries in a given field~$F$. Field~$F$ is supposed below to be \emph{finite}; more specifically, we will be dealing, in our examples in Section~\ref{s:e}, with the smallest Galois field $F=\mathbb F_2$ of two elements.

Let $R$ be one of the mentioned matrices. First, we extend it to a $6\times 6$ matrix $R_{ijk}$, where $1\le i<j<k\le 6$, as follows: $R_{ijk}$ is the direct sum of~$R$ put in the intersection of rows $i,j,k$ and columns of the same numbers, with the identity $3\times 3$ matrix (at the rest of rows/columns). For instance, if
\begin{equation*}
R = \begin{pmatrix} r_{11} & r_{12} & r_{13} \\ r_{21} & r_{22} & r_{23} \\ r_{31} & r_{32} & r_{33} \end{pmatrix},
\end{equation*}
then
\begin{equation*}
R_{145} = \begin{pmatrix} r_{11} & 0 & 0 & r_{12} & r_{13} & 0 \\ 
                          0 & 1 & 0 & 0 & 0 & 0 \\
                          0 & 0 & 1 & 0 & 0 & 0 \\
                          r_{21} & 0 & 0 & r_{22} & r_{23} & 0 \\ 
                          r_{31} & 0 & 0 & r_{32} & r_{33} & 0 \\
                          0 & 0 & 0 & 0 & 0 & 1
                          \end{pmatrix} .
\end{equation*}

The equation reads:
\begin{equation}\label{d}
R_{123}^{(1)} R_{145}^{(2)} R_{246}^{(3)} R_{356}^{(4)} = R_{356}^{(4)} R_{246}^{(3)} R_{145}^{(2)} R_{123}^{(1)} .
\end{equation}

From each matrix~$R$, a \emph{quantum} ``permutation-type'' matrix~$\mathcal R$ can be made, as explained in~\cite{Hietarinta}. Namely, consider the $\mathbb R$- or $\mathbb C$-linear space~$\mathcal V$, whose \emph{basis} consists of all elements of~$F$. In other words, $\mathcal V$ consists of formal $\mathbb R$- or $\mathbb C$-linear combinations of elements of~$F$.

The basis of the tensor product $\mathcal V^{\otimes 3} \stackrel{\mathrm{def}}{=} \mathcal V \otimes \mathcal V \otimes \mathcal V$ consists then of elements of $F^3$---that is, \emph{triples} of elements of~$F$. We prefer to think of these triples as organized in 3-\emph{rows}, on which matrix~$R$ can act from the \emph{right}; we follow thus the preferred way of dealing with vectors and matrices in the computer algebra system GAP~\cite{GAP}, in which our calculations were done.

$\mathcal R$ is, by definition, (the matrix of) the $\mathbb R$- or $\mathbb C$-linear operator sending each basis element~$\vec x \in F^3$ of~$\mathcal V^{\otimes 3}$ into another basis element, namely~$\vec x R$. Next, we consider the tensor product $\mathcal V^{\otimes 6}$ of \emph{six} copies~$\mathcal V_i$ of~$\mathcal V$, \ $i=1,\ldots,6$, and define $\mathcal R_{ijk}$ as the operator~$\mathcal R$ acting in ~$\mathcal V_i \otimes \mathcal V_j \otimes \mathcal V_k$, and tensor multiplied by the identity operator acting in the three remaining spaces (again like in~\cite{Hietarinta}). The quantum equation corresponding to (and following from)~\eqref{d} reads
\begin{equation}\label{q}
\mathcal R_{123}^{(1)} \mathcal R_{145}^{(2)} \mathcal R_{246}^{(3)} \mathcal R_{356}^{(4)} = \mathcal R_{356}^{(4)} \mathcal R_{246}^{(3)} \mathcal R_{145}^{(2)} \mathcal R_{123}^{(1)} .
\end{equation}

\section{Linear combinations of permutation-type $\mathcal R$-operators}\label{s:l}

Suppose that the following relation holds \emph{together} with~\eqref{d}: 
\begin{equation}\label{d2}
R_{123}^{(1)} R_{145}^{(2)} R_{246}^{(3)} S_{356}^{(4)} = S_{356}^{(4)} R_{246}^{(3)} R_{145}^{(2)} R_{123}^{(1)},
\end{equation}
where $S^{(4)}$ is another $R$-matrix,
\begin{equation}\label{4ne}
 S^{(4)}\ne R^{(4)}.
\end{equation}
Then, for the corresponding permutation-type quantum matrices,
\begin{equation}\label{q2}
\mathcal R_{123}^{(1)} \mathcal R_{145}^{(2)} \mathcal R_{246}^{(3)} \mathcal S_{356}^{(4)} = \mathcal S_{356}^{(4)} \mathcal R_{246}^{(3)} \mathcal R_{145}^{(2)} \mathcal R_{123}^{(1)}
\end{equation}
holds. A simple but remarkable fact is that, \emph{unlike} the corresponding direct-sum relations \eqref{d} and~\eqref{d2}, quantum relations \eqref{q} and~\eqref{q2}, if added together with arbitrary coefficients (from $\mathbb R$ or~$\mathbb C$), give again a quantum tetrahedron relation:
\begin{equation}\label{q+q2}
\mathcal R_{123}^{(1)} \mathcal R_{145}^{(2)} \mathcal R_{246}^{(3)} (\lambda\mathcal R_{356}^{(4)} + \mu\mathcal S_{356}^{(4)}) = (\lambda\mathcal R_{356}^{(4)} + \mu\mathcal S_{356}^{(4)}) \mathcal R_{246}^{(3)} \mathcal R_{145}^{(2)} \mathcal R_{123}^{(1)}
\end{equation}

Suppose now that, \emph{in addition} to \eqref{d} and~\eqref{d2}, the following two relations hold:
\begin{align}
R_{123}^{(1)} R_{145}^{(2)} S_{246}^{(3)} R_{356}^{(4)} = R_{356}^{(4)} S_{246}^{(3)} R_{145}^{(2)} R_{123}^{(1)}, \label{d3} \\
R_{123}^{(1)} R_{145}^{(2)} S_{246}^{(3)} S_{356}^{(4)} = S_{356}^{(4)} S_{246}^{(3)} R_{145}^{(2)} R_{123}^{(1)}, \label{d4}
\end{align}
where
\begin{equation}\label{3ne}
 S^{(3)}\ne R^{(3)}.
\end{equation}
Then, the following quantum tetrahedron relation holds:
\begin{equation}\label{q+q2+q3+q4}
 \begin{aligned}
\mathcal R_{123}^{(1)} & \mathcal R_{145}^{(2)} (\alpha\mathcal R_{246}^{(3)} + \beta\mathcal S_{246}^{(3)} ) (\lambda\mathcal R_{356}^{(4)} + \mu\mathcal S_{356}^{(4)}) \\
& = (\lambda\mathcal R_{356}^{(4)} + \mu\mathcal S_{356}^{(4)}) (\alpha\mathcal R_{246}^{(3)} + \beta\mathcal S_{246}^{(3)} ) \mathcal R_{145}^{(2)} \mathcal R_{123}^{(1)} .
 \end{aligned}
\end{equation}

The nontrivial part of two parameters $\alpha$ and~$\beta$ is their \emph{ratio}~$\frac{\beta}{\alpha}$, that is, \emph{one} parameter, belonging to the projective space $\mathbb {CP}^1 = \mathbb {C} \cup \{\infty\}$ or $\mathbb {RP}^1 = \mathbb {R} \cup \{\infty\}$; the similar remark holds for $\lambda$ and~$\mu$.

\section{Examples}\label{s:e}

Computer search for six-tuples
\begin{equation}\label{6t}
 \left( R^{(1)}\!,\: R^{(2)}\!,\: R^{(3)}\!,\: R^{(4)}\!,\: S^{(3)}\!,\: S^{(4)} \right)
\end{equation}
of matrices satisfying relations \eqref{d}, \eqref{d2}, \eqref{d3}, \eqref{d4}, with conditions \eqref{4ne} and~\eqref{3ne}, proves to be feasible on a personal computer, at least if we take the smallest field $F=\mathbb F_2$, and restrict ourself to the case where all matrices in~\eqref{6t} are \emph{invertible}, that is, belong to the group~$\mathrm{GL}(3,\mathbb F_2)$. The reason for the latter condition is just that $\mathrm{GL}(3,\mathbb F_2)$ has 168 elements, while there are considerably more, namely $2^9=512$, different $3\times 3$ matrices with entries in~$\mathbb F_2$.

There are 61535 solutions to equation~\eqref{d} with all~$R$'s belonging to~$\mathrm{GL}(3,\mathbb F_2)$; obtaining their list was the most time-consuming part of the work. Using that list, it was much easier to extract from it a list of all six-tuples~\eqref{6t}. Formally speaking, their number, namely 3828292, is quite impressive, but one must take into account that six-tuples obtained from each other by interchanges $R^{(3)}\leftrightarrow S^{(3)}$ and/or $R^{(4)}\leftrightarrow S^{(4)}$ give the same (families of) quantum matrices.

From these 3828292 six-tuples, of interest are \emph{nontrivial} ones: at least one of quantum matrices $\alpha\mathcal R^{(3)} + \beta\mathcal S^{(3)}$ and $\lambda\mathcal R^{(4)} + \mu\mathcal S^{(4)}$ must be \textbf{genuinely three-dimen\-sional}. The exact form of this condition that we adopt here is as follows: at least one of the matrices $R^{(3)}$, $S^{(3)}$, $R^{(4)}$ and $S^{(4)}$ must not be triangular or block triangular, even if we make a permutation of its rows and the same permutation of its columns (so, for instance, a matrix of the form $\begin{pmatrix} * & * & * \\ 0 & 1 & 0 \\ * & * & * \end{pmatrix}$ is \emph{not} suitable for this role). This simple condition ensures the interlinking of all three dimensions corresponding to this matrix, and hence to the quantum matrix, although weaker conditions making formal the idea of genuine three-dimen\-sionality may also be of interest.

The entries $0$ and~$1$ in the matrices below are understood as the elements of~$\mathbb F_2$. A statement like ``$\alpha \mathcal R^{(3)} + \beta \mathcal S^{(3)}$ gives a 14-vertex $\mathcal R$-matrix'' means that, for generic $\alpha$ and~$\beta$, this quantum matrix has 14 nonzero entries. If both $\alpha$ and~$\beta$ are positive, these entries are, clearly, also positive.

\begin{example}
 \begin{equation}
 \begin{aligned}
  R^{(1)} = \begin{pmatrix} 1 & 0 & 0 \\ 0 & 1 & 0 \\ 0 & 1 & 1 \end{pmatrix}, & \qquad
  R^{(2)} = \begin{pmatrix} 1 & 0 & 0 \\ 0 & 1 & 1 \\ 0 & 0 & 1 \end{pmatrix}, \\ 
  & R^{(3)} = \begin{pmatrix} 0 & 1 & 1 \\ 0 & 0 & 1 \\ 1 & 1 & 0 \end{pmatrix}, \qquad   
  R^{(4)} = \begin{pmatrix} 1 & 0 & 1 \\ 0 & 1 & 0 \\ 0 & 1 & 1 \end{pmatrix}, \qquad  
   \\
  & S^{(3)} = \begin{pmatrix} 1 & 0 & 0 \\ 0 & 1 & 0 \\ 0 & 0 & 1 \end{pmatrix}, \qquad  
  S^{(4)} = \begin{pmatrix} 1 & 1 & 1 \\ 0 & 1 & 0 \\ 0 & 1 & 1 \end{pmatrix} .  
 \end{aligned}
 \end{equation}
$\alpha \mathcal R^{(3)} + \beta \mathcal S^{(3)}$ gives a 14-vertex $\mathcal R$-matrix, while $\lambda \mathcal R^{(4)} + \mu \mathcal S^{(4)}$ gives a 12-vertex $\mathcal R$-matrix.
\end{example}

\begin{example}
 \begin{equation}
 \begin{aligned}
  R^{(1)} = \begin{pmatrix} 1 & 0 & 1 \\ 0 & 1 & 1 \\ 0 & 0 & 1 \end{pmatrix}, & \qquad
  R^{(2)} = \begin{pmatrix} 1 & 0 & 0 \\ 0 & 1 & 0 \\ 0 & 1 & 1 \end{pmatrix}, \\ 
  & R^{(3)} = \begin{pmatrix} 1 & 1 & 0 \\ 1 & 1 & 1 \\ 1 & 0 & 0 \end{pmatrix}, \qquad   
  R^{(4)} = \begin{pmatrix} 1 & 0 & 0 \\ 0 & 1 & 1 \\ 1 & 0 & 1 \end{pmatrix}, \qquad  
   \\
  & S^{(3)} = \begin{pmatrix} 0 & 0 & 1 \\ 1 & 1 & 1 \\ 1 & 0 & 0 \end{pmatrix}, \qquad  
  S^{(4)} = \begin{pmatrix} 1 & 0 & 0 \\ 1 & 1 & 1 \\ 1 & 0 & 1 \end{pmatrix} . 
 \end{aligned}
 \end{equation}
Both $\alpha \mathcal R^{(3)} + \beta \mathcal S^{(3)}$ and $\lambda \mathcal R^{(4)} + \mu \mathcal S^{(4)}$ give 12-vertex $\mathcal R$-matrices.
\end{example}

\begin{example}
 \begin{equation}
 \begin{aligned}
  R^{(1)} = \begin{pmatrix} 1 & 0 & 0 \\ 0 & 1 & 1 \\ 0 & 0 & 1 \end{pmatrix}, & \qquad
  R^{(2)} = \begin{pmatrix} 1 & 0 & 0 \\ 0 & 1 & 0 \\ 1 & 0 & 1 \end{pmatrix}, \\ 
  & R^{(3)} = \begin{pmatrix} 0 & 1 & 1 \\ 1 & 1 & 1 \\ 1 & 1 & 0 \end{pmatrix}, \qquad   
  R^{(4)} = \begin{pmatrix} 1 & 0 & 0 \\ 0 & 1 & 0 \\ 1 & 0 & 1 \end{pmatrix}, \qquad  
   \\
  & S^{(3)} = \begin{pmatrix} 1 & 0 & 0 \\ 1 & 1 & 1 \\ 1 & 1 & 0 \end{pmatrix}, \qquad  
  S^{(4)} = \begin{pmatrix} 1 & 0 & 0 \\ 1 & 1 & 0 \\ 1 & 0 & 1 \end{pmatrix} . 
 \end{aligned}
 \end{equation}
Both $\alpha \mathcal R^{(3)} + \beta \mathcal S^{(3)}$ and $\lambda \mathcal R^{(4)} + \mu \mathcal S^{(4)}$ give 12-vertex $\mathcal R$-matrices.
\end{example}

\begin{example}
 \begin{equation}
 \begin{aligned}
  R^{(1)} = \begin{pmatrix} 1 & 0 & 0 \\ 0 & 1 & 0 \\ 1 & 0 & 1 \end{pmatrix}, & \qquad
  R^{(2)} = \begin{pmatrix} 1 & 0 & 0 \\ 0 & 1 & 0 \\ 1 & 0 & 1 \end{pmatrix}, \\ 
  & R^{(3)} = \begin{pmatrix} 0 & 1 & 0 \\ 1 & 0 & 0 \\ 0 & 0 & 1 \end{pmatrix}, \qquad   
  R^{(4)} = \begin{pmatrix} 0 & 1 & 0 \\ 0 & 1 & 1 \\ 1 & 1 & 1 \end{pmatrix}, \qquad  
   \\
  & S^{(3)} = \begin{pmatrix} 1 & 1 & 0 \\ 0 & 1 & 0 \\ 0 & 0 & 1 \end{pmatrix}, \qquad  
  S^{(4)} = \begin{pmatrix} 0 & 1 & 1 \\ 1 & 0 & 1 \\ 1 & 1 & 1 \end{pmatrix} . 
 \end{aligned}
 \end{equation}
Both $\alpha \mathcal R^{(3)} + \beta \mathcal S^{(3)}$ and $\lambda \mathcal R^{(4)} + \mu \mathcal S^{(4)}$ give 14-vertex $\mathcal R$-matrices.
\end{example}

\section{One more way: modified tetrahedron}\label{s:m}

\subsection{Pairs of modified equations}\label{ss:pairs}

Suppose that, instead of~\eqref{d}, \emph{two} following ``modified'' equations hold:
\begin{align} 
R_{123}^{(1)} R_{145}^{(2)} R_{246}^{(3)} R_{356}^{(4)} & = Q_{356}^{(4)} R_{246}^{(3)} R_{145}^{(2)} R_{123}^{(1)} , \label{dd1} \\
R_{123}^{(1)} R_{145}^{(2)} R_{246}^{(3)} Q_{356}^{(4)} & = R_{356}^{(4)} R_{246}^{(3)} R_{145}^{(2)} R_{123}^{(1)} , \label{dd2} \\
R^{(4)} & \ne Q^{(4)}. \nonumber
\end{align}
Equations \eqref{dd1} and~\eqref{dd2} have, this way, the same triple $(R^{(1)},R^{(2)},R^{(3)})$, and turn into each other by an interchange $R^{(4)} \leftrightarrow Q^{(4)}$. We can write quantum equations 
\begin{align}
\mathcal R_{123}^{(1)} \mathcal R_{145}^{(2)} \mathcal R_{246}^{(3)} \mathcal R_{356}^{(4)} = \mathcal Q_{356}^{(4)} \mathcal R_{246}^{(3)} \mathcal R_{145}^{(2)} \mathcal R_{123}^{(1)} , \label{m1} \\
\mathcal R_{123}^{(1)} \mathcal R_{145}^{(2)} \mathcal R_{246}^{(3)} \mathcal Q_{356}^{(4)} = \mathcal R_{356}^{(4)} \mathcal R_{246}^{(3)} \mathcal R_{145}^{(2)} \mathcal R_{123}^{(1)} \label{m2}
\end{align}
corresponding to \eqref{dd1} and~\eqref{dd2} in the same way as we did before, and \emph{add them up}:
\begin{equation}\label{md}
\mathcal R_{123}^{(1)} \mathcal R_{145}^{(2)} \mathcal R_{246}^{(3)} \mathcal T_{356}^{(4)} = \mathcal T_{356}^{(4)} \mathcal R_{246}^{(3)} \mathcal R_{145}^{(2)} \mathcal R_{123}^{(1)} , 
\end{equation}
where
\begin{equation}\label{T}
\mathcal T^{(4)} = \mathcal R^{(4)} + \mathcal Q^{(4)} .
\end{equation}
The quantum operator~\eqref{T} will have \emph{integers} $0$, $1$ and~$2$ as its matrix entries.

\subsection{Examples}\label{ss:me}

Finding all pairs of relations~\eqref{dd1},~\eqref{dd2} is a considerably more difficult task for a personal computer than finding all relations~\eqref{d}. Nevertheless, not only examples of such pairs have been found, but striking examples have been discovered of the existence of \emph{two} or even \emph{four different} pairs $(R^{(4)}, Q^{(4)})$ for one given triple $( R^{(1)}, R^{(2)}, R^{(3)} )$. Hence, in such cases there exists, of course, a pair or a quadruple of different quantum~$\mathcal T$'s~\eqref{T} for one triple of quantum~$\mathcal R$'s and, moreover, any $\mathbb R$- or $\mathbb C$-linear combination of these~$\mathcal T$'s will again satisfy the tetrahedron together with the same triple of quantum~$\mathcal R$'s.

\begin{example}[One $\mathcal T$ for a triple of~$\mathcal R$'s]
\begin{equation}\label{t2}
R^{(1)} = \begin{pmatrix} 1 & 0 & 1 \\ 1 & 0 & 0 \\ 1 & 1 & 1 \end{pmatrix}, \qquad  
R^{(2)} = \begin{pmatrix} 1 & 0 & 0 \\ 0 & 1 & 1 \\ 0 & 0 & 1 \end{pmatrix}, \qquad  
R^{(3)} = \begin{pmatrix} 1 & 0 & 0 \\ 0 & 1 & 0 \\ 1 & 1 & 1 \end{pmatrix}.
\end{equation}

\begin{equation}
R^{(4)} = \begin{pmatrix} 1 & 0 & 0 \\ 0 & 1 & 0 \\ 1 & 0 & 1 \end{pmatrix}, \qquad  
Q^{(4)} = \begin{pmatrix} 1 & 0 & 0 \\ 0 & 1 & 0 \\ 1 & 1 & 1 \end{pmatrix}.
\end{equation}
\end{example}

\begin{example}[Two $\mathcal T$'s for a triple of~$\mathcal R$'s]
\begin{equation}\label{t4}
R^{(1)} = \begin{pmatrix} 1 & 0 & 0 \\ 0 & 1 & 0 \\ 1 & 0 & 1 \end{pmatrix}, \qquad 
R^{(2)} = \begin{pmatrix} 1 & 0 & 1 \\ 1 & 1 & 1 \\ 0 & 0 & 1 \end{pmatrix}, \qquad 
R^{(3)} = \begin{pmatrix} 1 & 0 & 1 \\ 1 & 1 & 0 \\ 0 & 0 & 1 \end{pmatrix}.
\end{equation}
       
Here are two pairs $(R^{(4)}, Q^{(4)})$ satisfying \eqref{dd1} and~\eqref{dd2} with the triple~\eqref{t4}.

First pair:       
\begin{equation}
R^{(4)} = \begin{pmatrix} 1 & 0 & 0 \\ 0 & 1 & 0 \\ 0 & 0 & 1 \end{pmatrix}, \qquad 
Q^{(4)} = \begin{pmatrix} 1 & 1 & 0 \\ 0 & 1 & 0 \\ 0 & 0 & 1 \end{pmatrix}.
\end{equation}
      
Second pair:       
\begin{equation}
R^{(4)} = \begin{pmatrix} 1 & 0 & 1 \\ 0 & 1 & 0 \\ 0 & 0 & 1 \end{pmatrix}, \qquad 
Q^{(4)} = \begin{pmatrix} 1 & 1 & 1 \\ 0 & 1 & 0 \\ 0 & 0 & 1 \end{pmatrix}.
\end{equation}

\end{example}

\begin{example}[One more example of two $\mathcal T$'s for a triple of~$\mathcal R$'s]

\begin{equation}\label{t4'}
R^{(1)} = \begin{pmatrix} 1 & 0 & 0 \\ 0 & 1 & 0 \\ 1 & 0 & 1 \end{pmatrix}, \qquad 
R^{(2)} = \begin{pmatrix} 1 & 0 & 1 \\ 0 & 1 & 0 \\ 0 & 0 & 1 \end{pmatrix}, \qquad 
R^{(3)} = \begin{pmatrix} 1 & 0 & 0 \\ 0 & 1 & 0 \\ 1 & 0 & 1 \end{pmatrix}. 
\end{equation}
       
Here are two pairs $(R^{(4)}, Q^{(4)})$ satisfying \eqref{dd1} and~\eqref{dd2} with the triple~\eqref{t4'}.

First pair:       
\begin{equation}      
R^{(4)} = \begin{pmatrix} 1 & 0 & 0 \\ 0 & 1 & 0 \\ 0 & 0 & 1 \end{pmatrix}, \qquad 
Q^{(4)} = \begin{pmatrix} 1 & 1 & 0 \\ 0 & 1 & 0 \\ 0 & 0 & 1 \end{pmatrix}.
\end{equation}

Second pair:       
\begin{equation}      
R^{(4)} = \begin{pmatrix} 1 & 1 & 0 \\ 0 & 1 & 0 \\ 0 & 1 & 1 \end{pmatrix}, \qquad 
Q^{(4)} = \begin{pmatrix} 1 & 0 & 0 \\ 0 & 1 & 0 \\ 0 & 1 & 1 \end{pmatrix}.
\end{equation}
\end{example}

\begin{example}[Four $\mathcal T$'s for a triple of~$\mathcal R$'s]

\begin{equation}\label{t8}
R^{(1)} = \begin{pmatrix} 1 & 0 & 0 \\ 0 & 1 & 0 \\ 1 & 0 & 1 \end{pmatrix}, \qquad
R^{(2)} = \begin{pmatrix} 1 & 0 & 1 \\ 0 & 1 & 0 \\ 0 & 0 & 1 \end{pmatrix}, \qquad 
R^{(3)} = \begin{pmatrix} 0 & 1 & 0 \\ 1 & 0 & 0 \\ 0 & 0 & 1 \end{pmatrix}.
\end{equation}

Here are four pairs $(R^{(4)}, Q^{(4)})$ satisfying \eqref{dd1} and~\eqref{dd2} with the triple~\eqref{t8}.

First pair:       
\begin{equation}
R^{(4)} = \begin{pmatrix} 1 & 0 & 0 \\ 0 & 1 & 0 \\ 0 & 0 & 1 \end{pmatrix}, \qquad 
Q^{(4)} = \begin{pmatrix} 1 & 1 & 0 \\ 0 & 1 & 0 \\ 0 & 0 & 1 \end{pmatrix},
\end{equation}

Second pair:
\begin{equation}
R^{(4)} = \begin{pmatrix} 1 & 1 & 0 \\ 0 & 1 & 0 \\ 0 & 1 & 1 \end{pmatrix}, \qquad 
Q^{(4)} = \begin{pmatrix} 1 & 0 & 0 \\ 0 & 1 & 0 \\ 0 & 1 & 1 \end{pmatrix},
\end{equation}

Third pair:      
\begin{equation}
R^{(4)} = \begin{pmatrix} 1 & 0 & 1 \\ 0 & 1 & 0 \\ 0 & 0 & 1 \end{pmatrix}, \qquad 
Q^{(4)} = \begin{pmatrix} 1 & 1 & 1 \\ 0 & 1 & 0 \\ 0 & 0 & 1 \end{pmatrix}.
\end{equation}

Fourth pair:      
\begin{equation}
R^{(4)} = \begin{pmatrix} 1 & 0 & 1 \\ 0 & 1 & 0 \\ 0 & 1 & 1 \end{pmatrix}, \qquad 
Q^{(4)} = \begin{pmatrix} 1 & 1 & 1 \\ 0 & 1 & 0 \\ 0 & 1 & 1 \end{pmatrix}.
\end{equation}

\end{example}

\section{Discussion}\label{s:hm}

An unusual feature of our solutions~\eqref{q+q2+q3+q4} is that $\alpha \mathcal R^{(3)} + \beta \mathcal S^{(3)}$ and $\lambda \mathcal R^{(4)} + \mu \mathcal S^{(4)}$ are \emph{different} families of $\mathcal R$-matrices. This may result in the appearance of an unusual algebra when studying integrable statistical physics model built from these $\mathcal R$-operators.

Concerning solutions~\eqref{md}, all $\mathcal R$-operators in them are, generally, different as well.

\end{document}